\documentclass[11pt]{article} 

\usepackage{amsthm}
\usepackage{amsmath}
\usepackage{amsfonts}
\usepackage{amssymb}
\usepackage{amscd}
\usepackage{graphicx}

\usepackage{tikz-cd}
\usetikzlibrary{matrix,arrows,decorations.pathmorphing}

\newtheorem{theorem}{Theorem}[section]
\newtheorem{lemma}[theorem]{Lemma}

\newtheorem{definition}[theorem]{Definition}
\newtheorem{example}[theorem]{Example}

\newtheorem{remark}[theorem]{Remark}

\newtheorem{corollary}[theorem]{Corollary}
\usepackage{geometry} 
\usepackage{amssymb}
\usepackage{graphicx} 

\usepackage{float} 
\usepackage{wrapfig} 

\usepackage{lipsum} 

\graphicspath{{./Pictures/}} 
\usepackage{tikz}

  \usepackage{txfonts}
\usepackage[arrow, matrix, curve]{xy}
\usepackage{array}
\usepackage{setspace}
\usepackage{mathrsfs}
\newcommand{\IZ}{\mathbb{Z}}

\newcommand{\IP}{\mathbb{P}}

\newtheorem*{thm}{Theorem}

\numberwithin{table}{section}
\begin{document}

\title{The Locus of Curves with $D_n$-Symmetry inside $\mathfrak M_g$ }
\date{}
\author{Binru Li, Sascha Weigl}
\maketitle

\begin{abstract}
\noindent
The aim of this paper is to determine the irreducible components of  $\mathfrak{M}_{g}(D_n)$, the locus  inside $\mathfrak M_g$ of the curves admitting an effective action by the dihedral group $D_n$. 
This is done by classifying pairs $(H,H')$ of distinct subgroups of the mapping class group $Map_g$, such that both $H$ and $H'$ are isomorphic to $D_n$ and the fixed point locus of $H$ inside the Teichm\"uller space $\mathcal T_g$ is contained in the fixed point locus of $H'$.

\end{abstract}

\section{Introduction}
Given a finite group $H$, denote by $\mathfrak{M}_g(H)$ the locus inside $\mathfrak{M}_g$ (the coarse moduli space of curves of genus $g\geq 2$) of the curves admitting an effective action by the group $H$. A good approach to understanding the irreducible components of $\mathfrak{M}_g(H)$ is to view $\mathfrak{M}_g$ as the quotient of the Teichm\"uller space $\mathcal T_g$ by the natural action of the mapping class group $Map_g$: 
$$\pi: \mathcal T_g \to \mathcal T_g/Map_g = \mathfrak M_g .$$ 
 Observe that $$\mathfrak{M}_{g}(H) = \bigcup \limits_{[\rho]} \mathfrak M_{g,~\rho}(H),$$
where $\rho:H\hookrightarrow Map_g$ is an injective homomorphism, $\mathfrak{M}_{g,\ \rho}(H)$ is the image of the fixed locus of $\rho(H)$ under the natural projection $\pi$ and $\rho \sim \rho'$ iff they are equivalent by the equivalence relation generated by the automorphisms of $H$ and the conjugations by $Map_g$. We call this equivalence class an {\em unmarked topological type} (cf. \cite{CLP2}, section 2). Since each $\mathfrak M_{g,~\rho}(H)$ is an irreducible (Zariski) closed subset of $\mathfrak{M}_g$ (cf. \cite{CLP2}, Theorem 2.3), in order to determine the irreducible components of $\mathfrak{M}_{g}(H)$, it suffices to determine the maximal loci of the form $\mathfrak M_{g,~\rho}(H)$, i.e. to figure out when one locus contains another. \\
The case where $H$ is a cyclic group was investigated in \cite{Cor} and \cite{Cat1}. In \cite{CLP2} the authors have defined a new homological invariant which allows them to tell when two homomorphism $\rho$ and $\rho'$ are not equivalent; for the case of $H=D_n$, the dihedral group, they also found one representative for each unmarked topological type.\\
In this paper, we focus on the case $H=D_n$, and solve the following problem: for which $\rho$ and $\rho'$, does $\mathfrak M_{g,~\rho'}(D_n)$ contain $\mathfrak M_{g,~\rho}(D_n)$? Hence we determine the loci $\mathfrak{M}_{g,\ \rho}(D_n)$ which are not maximal whence the irreducible decomposition of $\mathfrak{M}_g(D_n)$. The above problem is equivalent to the classification of subgroups $H,H'$ of $Map_g$ $(g\geq 2)$, where $H$ and $H'$ satisfy the following condition: $$(*)~H,H'  \simeq D_n, H\neq H' ~\mbox{and}~Fix(H)\subset Fix(H') .$$

\noindent
For any finite subgroup $H\subset Map_g$, set $\delta_H:=$dim$Fix(H)$ and let $G:=G(H):=\bigcap_{C\in Fix(H)}Aut(C)$ ($Fix(H)$ corresponds to the complex structures for which the action of $H$ is holomorphic, whereas $G(H)$ is the common automorphism group of all the curves in $Fix(H)$). 
 If $H=G(H)$ we call $H$ {\em full}.\\
It is easy to see that condition $(*)$ is equivalent to the condition\\
$(**)$ $H$ is isomorphic to $D_n$ and not full, $G(H)$ has a subgroup $H'$ which is isomorphic to $D_n$ and different from $H$.  \\

For any curve $C\in Fix(G)$, we have a Galois cover $p: C\rightarrow C/G=:C'$ which is branched in $r$ ($r$ can be zero) points $P_1,...,P_r$ on $C'$ with branching indices $m_1,...,m_r$. By Theroem \ref{MSSV}, in our case $C'$ is always $\IP^1$. The cover map $p$ is determined by a surjective homomorphism $f$ from the orbifold fundamental group $T(m_1,...,m_r):=<\gamma _1,...,\gamma _r |\Pi \gamma_i=1,\gamma_i^{m_i}=1>$ to $G$ (cf. \cite{Cat2}, section 5). The vector $v:=(f(\gamma_1),...,f(\gamma_r))$ is called the \emph{Hurwitz vector} associated to $f$ (See section 5 for more details). Then two Hurwitz vectors $v$ and $v'$ determine the same topological type if and only if they are equivalent for the equivalence relation generated by the action of $Aut(G)$ and by sequences of braid moves. (See Definition \ref{equivalent HV}).\\

Our main result is the following:

\begin{thm}
Let $H,H'$ be subgroups of $Map_g$, satisfying condition $(*)$. Then $G(H) \simeq D_n \times \IZ/2$ and $H$ corresponds to $D_n \times \{0\}$. The group $H'$ and the topological action of the group $G(H)$ (i.e. its Hurwitz vector) are as listed in the tables of section 2.
\end{thm}

\noindent
The structure of this paper is as follows:\bigskip

\noindent
In section 2 we present our results through tables.\\
In section 3 we quote a Theorem from \cite{MSSV} (cf. Theorem \ref{MSSV}), which contains the possible cases (which we call \emph{cover type}) where $H \subsetneq G \subset Map_g$ and $\delta_G=\delta_H$. From this Theorem, using the Riemann-Hurwitz formula, we obtain pairs of dimensions $(\delta_H,\delta_{H'})$, which can occur under condition (**).
 We will also see that $C/G\simeq \mathbb{P}^1$ and $[G:H]=2$ except for one case.\\
In section 4 we will understand group theoretically which cases of $H$ and $G$ can happen under condition (**). This is done by classifying the index 2 subgroups of $G$, where $G$ is a finite group containing two distinct index 2 subgroups which are isomorphic to $D_n$. The cases there are called the \textit{group types}. \\
In section 5 we classify the equivalence classes of Hurwitz vectors of the map $C\rightarrow C/G\simeq \mathbb{P}^1$ for each cover type and group type, by giving one representative vector for each equivalence class.

\section{Results}
We present our results through tables. There will be one table for each normal form of Hurwitz vectors for the covering $C \to C/G$, obtained in section 5. For the reader's convenience we present a short list of notation:\\


\begin{tabular}{cl}

$\nu$ & Hurwitz vector for the covering $C \to C/G$ \\

$v_{G/H'}$ &  Hurwitz vector for the double covering $C/H' \to C/G = \IP^1$ \\

$g_{C/H'}$& Genus of $C/H'$ \\

$\delta_{H'}$& Dimension of $Fix(H')$\\

$v_{H'}$& Hurwitz vector for the covering $C \to C/H'$

\end{tabular}

\vspace{0,5cm}

\noindent
We will use the following subgroups of $D_n \times \IZ/2$, where  $D_n =<x,y~|~x^n=y^2=1, yxy^{-1} = x^{-1}>$ and $e$ denotes the neutral element of $D_n$.\\

\begin{center}
\begin{tabular}{c|c}
Subgroup&Generators\\
\hline
$K$ & $(x,0)$\\

$H_{1,1}$ & $K,(e,1)$\\

$H_{1,2}$ & $K,(y,1)$\\

$H_{1,3}$& $(x^2,0),(y,0),(e,1)$\\

$H_{1,4}$&$(x^2,0),(y,0),(x,1)$\\

$H_{1,5}$&$(x^2,0),(yx,0),(e,1)$\\

$H_{1,6}$&$(x^2,0),(yx,0),(x,1)$\\

\end{tabular}
\end{center}

\bigskip
\noindent
For compactness, we make the following conventions:\\
 Whenever the groups $H_{1,4},H_{1,6},H_{1,3},H_{1,5}$ occur, we assume that $n = 2m$, in the last 2 cases we additionally assume $m$ to be odd. If $H_{1,1}$ appears we are in the case $n = 2$. We identify the groups $H_{1,3}$ and $H_{1,5}$ with $D_n$ by sending their respective generators in the given order to $x^{m+1},y,x^m$.\\
The cover types are those which appear in Theorem \ref{MSSV}.

\begin{theorem}\label{main result}
Let $H,H'$ be subgroups of $Mapg_g$, satisfying condition $(*)$. Then $G(H) \simeq D_n \times \IZ/2$, $H$ corresponds to $D_n \times \{0\}$. The group $H'$ and the topological action of the group $G(H)$ (i.e. its Hurwitz vector) are as listed in the following tables.
\end{theorem}
We obtain immediately the following corollary:
\begin{corollary}
The locus $\mathfrak{M}_{g,\ \rho}(D_n)$ is maximal iff its topological type $[\rho]$ is different from those which are determined by $C\rightarrow C/H$ in the following tables.
\end{corollary}
\begin{remark}
Given a cover $C\rightarrow C/H$, the data consisting of $g_{C/H}$ and the branching indices are called the signature of the cover. In $\cite{BCGG}$, section 3 the authors computed the signatures for the possible non-maximal loci of the form $\mathfrak{M}_{g,\ \rho}(D_n)$, which is a corollary of our result.
\end{remark}

\vspace{1cm} 
 \begin{center} \begin{large}Cover type I)\\
  ($\delta_H =3$, $g_{C/H}=2$, $C\rightarrow C/H$ is unramified)\end{large} \end{center}
  \[v=(((y,1),(y,1),(yx,1),(yx,1),(e,1),(e,1))\]

\begin{center}
\begin{tabular}{|c|c|c|c|c|}
\hline
$H'$&$v_{G/H'}$&$g_{C/H'}$&$\delta_{H'}$&$v_{H'}$\\
\hline
$H_{1,2}$&(0,0,0,0,1,1)&0&5&(y,y,y,y,yx,yx,yx,yx)\\
\hline
$H_{1,3}$&(0,0,1,1,0,0)&0&5&$(yx^m,yx^{m-2},yx^m,yx^{m+2},x^m,x^m,x^m,x^m)$\\
\hline
$H_{1,4}$&(1,1,0,0,1,1)&0&5&(yx,yx,yx,yx,y,y,y,y)\\
\hline
$H_{1,5}$&(1,1,0,0,0,0)&0&5&$(yx^m,yx^m,yx^m,yx^m,x^m,x^m,x^m,x^m)$\\
\hline
$H_{1,6}$&(0,0,1,1,1,1)&1&4&(e,yx;y,y,y,y)\\
\hline
\end{tabular}
\end{center}





\[v = ((y,1),(yx^{m},1),(yx,1),(yx,1),(x^{m},1),(e,1)),\ n=2m\]

\begin{center}
\begin{tabular}{|c|c|c|c|c|}
\hline
$H'$&$v_{G/H'}$&$g_{C/H'}$&$\delta_{H'}$&$v_{H'}$\\
\hline
$H_{1,2}$&$(0,0,0,0,1,1)$&0&5&$(y,y,yx^m,x^my,yx,yx,yx,yx)$\\
\hline
$H_{1,3}$&$(0,1,1,1,1,0)$&1&4&$(x^{m+1},x^{m-1};x^m,x^m,yx^m,yx^m)$\\
\hline
$H_{1,4}$($m$ odd)&$(1,0,0,0,0,1)$&0&5&$(yx^m,x^my,yx,yx^3,yx,xy,x^m,x^m)$\\
$H_{1,4}$($m$ even)&$(1,1,0,0,1,1)$&1&4&$(x^m,x^my;yx,yx,yx,yx)$\\
\hline
$H_{1,5}$&$(1,0,0,0,1,0)$&0&5&$(yx^{\frac{m^2-1}{2}},yx^{\frac{m^2-1}{2}},yx^m,yx^m,yx^m,yx^m,x^m,x^m)$\\
\hline
$H_{1,6}$($m$ odd)&(0,1,1,1,0,1)&1&4&$(x^{m+1},x^{m-1};x^m,x^m,y,y)$\\
$H_{1,6}$($m$ even)&(0,0,1,1,1,1)&1&4&$(e,x^{m-1}y;y,y,x^my,yx^m)$\\
\hline

\end{tabular}
\end{center}

\newpage

\[v = ((y,1),(yx^{m},1),(yx^2,1),(yx^2,1),(x^{m},1),(e,1)),\ n=2m,\ m ~\mbox{odd.}\]

\begin{center}
\begin{tabular}{|c|c|c|c|c|}
\hline
$H'$&$v_{G/H'}$&$g_{C/H'}$&$\delta_{H'}$&$v_{H'}$\\
\hline
$H_{1,2}$&$(0,0,0,0,1,1)$&0&5&$(y,y,yx^m,yx^m,yx^2,yx^2,yx^2,yx^2)$\\
\hline
$H_{1,3}$&$(0,1,0,0,1,0)$&0&5&$(x^m,x^m,yx^m,yx^{-1},yx^{-3},yx^{-1},yx)$\\
\hline
$H_{1,4}$&$(0,0,0,0,1,1)$&0&5&$(y,y,yx^m,yx^m,yx^2,yx^2,yx^2,yx^2)$\\
\hline
$H_{1,5}$&$(1,0,1,1,1,0)$&1&4&$(x^2,x^{-2};x^m,x^m,x^my,x^my)$\\
\hline
$H_{1,6}$&(0,1,0,0,0,1)&0&5&$(y,y,x^2y,x^6y,x^2y,yx^2,x^m,x^m)$\\
\hline
\end{tabular}
\end{center}

For $n=2$ we have two extra cases:
\[ v = ((y,1),(y,1),(x,1),(x,1),(e,1),(e,1)) \]

\begin{center}
\begin{tabular}{|c|c|c|c|c|}
\hline
$H'$&$v_{G/H'}$&$g_{C/H'}$&$\delta_{H'}$&$v_{H'}$\\
\hline
$H_{1,1}$&$(1,1,0,0,0,0)$&0&5&$(x,x,x,x,y,y,y,y)$\\
\hline
$H_{1,2}$&$(0,0,1,1,1,1)$&1&4&$(e,x;y,y,y,y)$\\
\hline
$H_{1,3}$&$(0,0,1,1,0,0)$&0&5&$(yx,yx,yx,yx,x,x,x,x)$\\
\hline
$H_{1,4}$&$(1,1,0,0,1,1)$&1&4&$(e,y;x,x,x,x)$\\
\hline
$H_{1,5}$&$(1,1,1,1,0,0)$&1&4&$(e,yx;x,x,x,x)$\\
\hline
$H_{1,6}$&(0,0,0,0,1,1)&0&5&$(y,y,y,y,x,x,x,x)$\\
\hline
\end{tabular}
\end{center}

\[v = ((y,1),(yx,1),(x,1),(x,1),(x,1),(e,1))\]

\begin{center}
\begin{tabular}{|c|c|c|c|c|}
\hline
$H'$&$v_{G/H'}$&$g_{C/H'}$&$\delta_{H'}$&$v_{H'}$\\
\hline
$H_{1,1}$&$(1,1,0,0,0,0)$&0&5&$(yx,yx,yx,yx,yx,yx,y,y)$\\
\hline
$H_{1,2}$&$(0,0,1,1,1,1)$&1&4&$(e,e;y,y,yx,yx)$\\
\hline
$H_{1,3}$&$(0,1,1,1,1,0)$&1&4&$(y,y;x,x,yx,yx)$\\
\hline
$H_{1,4}$&$(1,0,0,0,0,1)$&0&5&$(yx,yx,x,x,x,x,x,x,x,x)$\\
\hline
$H_{1,5}$&$(1,0,1,1,1,0)$&1&4&$(e,y;x,x,x,x)$\\
\hline
$H_{1,6}$&(0,1,0,0,0,1)&0&5&$(y,y,x,x,x,x,x,x,x,x)$\\
\hline
\end{tabular}
\end{center}
\vspace{0.1cm}
\newpage
\begin{center} \begin{large}Cover type II) \\
($\delta_H=2$, $g_{C/H}=1$)\end{large} \end{center}
$(1)$ $c_5=2$.
$$v=((y,1),(yx,1),(yx,1),(e,1),(y,0)),\ v_H=(x,x^{-1};y,y).$$
\begin{center}
\begin{tabular}{|c|c|c|c|c|}
\hline
$H'$&$v_{G/H'}$&$g_{C/H'}$&$\delta_{H'}$&$v_{H'}$\\
\hline
$H_{1,2}$&$(0,0,0,1,1)$&0&3&$(y,y,yx,yx,yx,yx)$\\
\hline
$H_{1,3}$&$(0,1,1,0,0)$&0&3&$(yx^m,yx^{-1},x^m,x^m,y,yx^{m+1})$\\
\hline
$H_{1,4}$&$(1,0,0,0,1)$&0&3&$(yx,yx,yx,yx,y,y)$\\
\hline
$H_{1,5}$&$(1,0,0,0,1)$&0&3&$(yx^m,yx,yx^m,yx^{-1},x^m,x^m)$\\
\hline
$H_{1,6}$&$(0,1,1,1,1)$&1&2&$(e,yx;y,y)$\\ 
\hline
\end{tabular}
\end{center}

$(2)$ $c_5>2$.
$$v=((y,1),(yx^{-1},1),(e,1),(e,1),(x,0)),\ c_5=n,\ v_H=(x^{-1},y;x,x).$$
\begin{center}
\begin{tabular}{|c|c|c|c|c|}
\hline
$H'$&$v_{G/H'}$&$g_{C/H'}$&$\delta_{H'}$&$v_{H'}$\\
\hline
$H_{1,2}$&$(0,0,1,1,0)$&0&3&$(y,y,yx^{-1},yx^{-3},x,x)$\\
\hline
$H_{1,3}$&$(0,1,0,0,1)$&0&4&$(yx^m,yx^{-1},x^m,x^m,x^m,x^m,x^{m+1})$\\
\hline
$H_{1,4}$&$(1,0,1,1,1)$&1&3&$(y,y;x^2,yx^3,yx)$\\
\hline
$H_{1,5}$&$(1,0,0,0,1)$&0&4&$(yx^{-1},yx^{m-2},x^m,x^m,x^m,x^m,x^{m+1})$\\
\hline
$H_{1,6}$&$(0,1,1,1,1)$&1&3&$(yx^{-1},yx^{-1};x^2,yx^2,y)$\\ 
\hline
\end{tabular}
\end{center}

$$v=((y,1),(yx^{m-1},1),(x^m,1),(e,1),(x,0)),\ n=2m,\ v_H=(x^{m-1},yx^m;x,x).$$
\begin{center}
\begin{tabular}{|c|c|c|c|c|}
\hline
$H'$&$v_{G/H'}$&$g_{C/H'}$&$\delta_{H'}$&$v_{H'}$\\
\hline
$H_{1,2}$&$(0,0,1,1,0)$&0&3&$(y,y,yx^{m-1},yx^{m-3},x,x)$\\
\hline
$H_{1,3}$&$(0,0,1,0,1)$&0&4&$(yx^m,yx^{m+2},yx,yx,x^m,x^m,x^{-2})$\\
\hline
$H_{1,4}\ (m\ odd)$ &$(1,1,0,1,1)$&1&3&$(x^{m-1},y;x^2,x^m,x^m)$\\
$H_{1,4}\ (m\ even)$ &$(1,0,1,1,1)$&1&3&$(yx^{m},y;x^2,yx^{m+3},yx^{m+1})$\\
\hline
$H_{1,5}$&$(1,1,1,0,1)$&1&3&$(x^{m+1},y;x^{-2},x^m,x^m)$\\
\hline
$H_{1,6}\ (m\ odd)$&$(0,0,0,1,1)$&0&4&$(y,yx^{-2},yx^{m-1},yx^{m-1},x^m,x^m,x^2)$\\ 
$H_{1,6}\ (m\ even)$&$(0,1,1,1,1)$&1&3&$(yx^{-1},yx^{m-1};x^2,yx^2,y)$\\ 
\hline
\end{tabular}
\end{center}
\newpage
$$v=((y,1),(yx^{m-2},1),(x^m,1),(e,1),(x^2,0)),\ n=2m,\ m\ \mbox{odd},\ c_5=m,$$
 $$v_H=(x^{m-2},yx^m;x^2,x^2),$$
\begin{center}
\begin{tabular}{|c|c|c|c|c|}
\hline
$H'$&$v_{G/H'}$&$g_{C/H'}$&$\delta_{H'}$&$v_{H'}$\\
\hline
$H_{1,2}$&$(0,0,1,1,0)$&0&3&$(y,y,yx^{m-2},yx^{m-6},x^2,x^2)$\\
\hline
$H_{1,3}$&$(0,1,1,0,0)$&0&3&$(yx^m,yx^{m-4},x^m,x^m,x^2,x^2)$\\
\hline
$H_{1,4}$&$(1,0,0,1,0)$&0&3&$(yx^{m-2},yx^{m-6},x^m,x^m,x^2,x^2)$\\
\hline
$H_{1,5}$&$(1,0,1,0,0)$&0&3&$(y,yx^{-4},x^m,x^m,x^2,x^2)$\\
\hline
$H_{1,6}$&$(0,1,0,1,0)$&0&3&$(y,yx^{-4},x^m,x^m,x^2,x^2)$\\ 
\hline
\end{tabular}
\end{center}

For $n=2$ we have one extra case.\\
$$v=((yx,1),(x,1),(e,1),(e,1),(y,0)),\ v_H=(y,yx;y,y).$$
\begin{center}
\begin{tabular}{|c|c|c|c|c|}
\hline
$H'$&$v_{G/H'}$&$g_{C/H'}$&$\delta_{H'}$&$v_{H'}$\\
\hline
$H_{1,1}$&$(1,0,0,0,1)$&0&3&$(x,x,y,y,y,y)$\\
\hline
$H_{1,2}$&$(0,1,1,1,1)$&1&2&$(x,x;yx,yx)$\\
\hline
$H_{1,3}$&$(1,1,0,0,0)$&0&3&$(x,x,x,x,y,y)$\\
\hline
$H_{1,4}$&$(0,0,1,1,0)$&0&3&$(yx,yx,x,x,y,y)$\\
\hline
$H_{1,5}$&$(0,1,0,0,1)$&0&3&$(yx,yx,x,x,x,x)$\\
\hline
$H_{1,6}$&$(1,0,1,1,1)$&1&2&$(yx,yx;x,x)$\\ 
\hline
\end{tabular}
\end{center}

\begin{center} \begin{large}Cover type III-a)\\
 ($\delta_H =1$, $g_{C/H}=1$)\end{large} \end{center}
$$v=((y,1),(yx^{-1},1),(e,1),(x,1)),\ 2d_4=n=2m,\ v_H=(x^{-1},y;x^2).$$
\begin{center}
\begin{tabular}{|c|c|c|c|c|}
\hline
$H'$&$v_{G/H'}$&$g_{C/H'}$&$\delta_{H'}$&$v_{H'}$\\
\hline
$H_{1,2}$&$(0,0,1,1)$&0&2&$(y,yx^2,yx^{-1},yx^{-1},x^2)$\\
\hline
$H_{1,3}$&$(0,1,0,1)$&0&2&$(yx^m,yx^{-1},x^m,x^m,x^{m+1})$\\
\hline
$H_{1,4}$&$(1,0,1,0)$&0&1&$(yx^{-1},yx^{-3},x,x)$\\
\hline
$H_{1,5}$&$(1,0,0,1)$&0&2&$(yx^{-1},yx^{m-2},x^m,x^m,x^{m+1})$\\
\hline
$H_{1,6}$&$(0,1,1,0)$&0&1&$(x,x,y,yx^{-2})$\\ 
\hline
\end{tabular}
\end{center}
\newpage
\begin{center} \begin{large}Cover type III-b) \\
($\delta_H=1$, $g_{C/H}=0$)\end{large} \end{center}
$$v=((yx,1),(e,1),(y,0),(x,0)),\ c_4=n=2m,\ v_H=(y,yx^{-2},x,x).$$
\begin{center}
\begin{tabular}{|c|c|c|c|c|}
\hline
$H'$&$v_{G/H'}$&$g_{C/H'}$&$\delta_{H'}$&$v_{H'}$\\
\hline
$H_{1,2}$&$(0,1,1,0)$&0&1&$(yx,yx^{-1},x,x)$\\
\hline
$H_{1,3}$&$(1,0,0,1)$&0&2&$(x^m,x^m,y,yx^{m-1},x^{m+1})$\\
\hline
$H_{1,4}$&$(0,1,0,1)$&0&2&$(yx,yx^{-1},y,y,x^2)$\\
\hline
$H_{1,5}$&$(0,0,1,1)$&0&2&$(yx^m,yx^{-1},x^m,x^m,x^{m+1})$\\
\hline
$H_{1,6}$&$(1,1,1,1)$&1&1&$(yx,x;x^2)$\\ 
\hline
\end{tabular}
\end{center}
\bigskip
\section{A rough classification}
In this section we determine the possible pairs of dimensions $(\delta_H, \delta_{H'})$, for distinct subgroups $H$ and $H'$ of $Map_g$ which satisfy condition (**). \\
Given $C\in Fix(H)$, assume that $C\rightarrow C/H$ is a cover branched on $r$ points. We have that $\delta_H=3g_{G/H}-3+r$ (cf. \cite{CLP2}, Theorem 2.3).\\
The case $\delta_H=\delta_{H'}$ was done in Corollary 7.2 of [CLP2]. We only consider the case $\delta_H<\delta_{H'}$.\\
We recall Lemma 4.1 of [MSSV]:

\begin{theorem}\label{MSSV} (MSSV)\\
Let $H\subsetneqq G$ be two (finite) subgroups of $Map_g$, $\delta_H=\delta_G=:\delta$. Then one of the following holds:\\
$I)$ $\delta_{H}=3$, $[G$:$H]=2$, $C\rightarrow C/G$ is a covering of $\mathbb{P}^1$ branched on 6 points $P_1,\dots,P_6$, and with branching indices all equal to 2. Moreover the subgroup H corresponds to the unique genus two double cover of $\mathbb{P}^1$ branched on the 6 points.\\
$II)$ $\delta_H=2$, $[G$:$H]=2$, and $C\rightarrow C/G$ is a covering of $\mathbb{P}^1$ branched on five points, $P_1,\dots,P_5$, with branching indices $2,2,2,2,c_5$. Moreover the subgroup H corresponds to a double cover of $\mathbb{P}^1$ branched on the 4 points $P_1,\dots, P_4$ with branching index 2.\\
$III)$ $\delta_{H}=1$, there are 3 possibilities:\\
$III-a)$ H has index 2 in G, and $C\rightarrow C/G$ is a covering of $\mathbb{P}^1$ branched on 4 points, $P_1,\dots,P_4$, with branching indices $2,2,2,2d_4,$ where $d_4>1$. Moreover the subgroup H corresponds to the unique genus one double cover of $\mathbb{P}^1$ branched on the 4 points $P_1,\dots, P_4$.\\
$III-b)$ H has index 2 in G, and $C\rightarrow C/G$ is a covering of $\mathbb{P}^1$ branched on 4 points, $P_1,\dots,  P_4$, with branching indices $2,2,c_3,c_4$, where $c_3\leq c_4$ and $c_4>2$. Moreover the subgroup $H$ corresponds to a genus zero double cover of  $\mathbb{P}^1$ branched on two points with branching index 2.\\
$III-c)$ $H$ is normal in $G$, $G/H\cong (\mathbb{Z}/2)^2$, moreover $C\rightarrow C/G$ is a covering of $\mathbb{P}^1$ branched on 4 points $P_1,\dots, P_4$, with branching indices $2,2,2,c_4$, where $c_4>2$. Moreover the subgroup $H$ corresponds to the unique genus zero cover of $\mathbb{P}^1$ with group $(\mathbb{Z}/2)^2$ branched on the 3 points $P_1,P_2,P_3$ with branching index 2.
\end{theorem}
We call the cases in Theorem \ref{MSSV} the \emph{cover type} (of $H$ and $G$).\\
Since we have condition $(**)$, which implies $\delta_G = \delta_H$, we can apply Theorem \ref{MSSV}. Moreover we apply the Riemann-Hurwitz formula to each cover type to find the possible pairs $(\delta_H,\delta_{H'})$.

\begin{corollary}\label{MSSV Cor}
Assume $(**)$ and moreover $\delta_H<\delta_{H'}$. Then the following pairs of dimensions $(\delta_H,\delta_{H'})$ can occur:\\

\noindent
$I)$ $(3,4)$, $(3,5)$.\\
$II)$ $(2,3)$, $(2,4)$.\\
$III-a)$ $(1,2)$.\\
$III-b)$ $(1,2),(1,3)$.\\
$III-c)$ None.
\begin{proof}
$I)$ $\delta_H=3$.\\
 By the Riemann-Hurwitz formula,
$$2g(C)-2=|G|(-2+6 \cdot \frac{1}{2})=|H'|(2(g_{C/H'}-1)+k/2)$$
 where $k$ is the number of  branching points of $C\rightarrow C/H'$.\\
It is easy to see that $(g_{C/H'},k)=(2,0),(1,4)$ or $(0,8)$, corresponding to $\delta_{H'}=3,4,5$. Since we require $\delta_{H}<\delta_{H'}$, the possible pairs are (3,4) and (3,5).\\
$II)$ $\delta_H=2$.\\
 In this case $C/H'\rightarrow \mathbb{P}^1$ is a double covering branched on at most 5 points. Using Riemann-Hurwitz, there are two cases:\\
(i)  $g_{C/H'}=0$ and $C/H' \rightarrow \mathbb{P}^1$ is branched on 2 of the 5 points with branching indices 2,2. \\
If $c_5=2$ or $P_5$ is not a branching point, we have $\delta_{H'}=3$;\\
Otherwise $c_5$ is even and bigger than 2 and $P_5$ is a branching point, we get $\delta_{H'}=4$. \\
(ii) $g_{C/H'}=1$ and $C/H' \rightarrow \mathbb{P}^1$ is branched on 4 of the 5 points with branching indices 2,2,2,2.\\
The only possible case in which $\delta_{H'}>2$ is that $c_5$ is even and bigger than 2 and $P_5$ is one of the branching points. In this case $\delta_{H'}=3$.\\
$III)$ $\delta_H=1$.\\
$III-a)$
Similar to case $II)$, one gets $g_{C/H'}=0$, and $C/H'\rightarrow \mathbb{P}^1$ is a double cover with one of the branching points $P_4$ and $\delta_{H'}=2$.\\
$III-b)$
$i)$ If $c_3=2$, the only possibility is $c_4$ even, $g_{C/H'}=0$ and $C/H'\rightarrow \mathbb{P}^1$ is a double cover with one of the branching points $P_4$, here $\delta_{H'}=2$.\\
$ii)\  c_3>2$, there are three possibilities:\\
$\alpha)$ $c_3$ or $c_4$ is even, one and only one point of $P_3$,$P_4$ is a branching point. This case is similar to $III-b)-i)$, $\delta_{H'}=2$.\\
$\beta)$ Both $c_3$ and $c_4$ are even, $g_{C/H'}=0$, and $C/H'\rightarrow \mathbb{P}^1$ is a double cover branching on $P_3$,$P_4$. We have $\delta_{H'}=3$.\\
$\gamma)$ Both $c_3$ and $c_4$ are even, $g_{C/H'}=1$, and $C/H'\rightarrow \mathbb{P}^1$ is a double cover branching on 4 points $P_1,\dots, P_4$. We have $\delta_{H'}=2$.\\
$III-c)$ We will give the proof in section \ref{admissible f's}, Lemma \ref{group type 3}.
\end{proof}
\end{corollary}


Remark: Cor. \ref{MSSV Cor} is valid for any $H,H'$ with  the same index in $G$ except for the case $III-c)$.
\section{Index 2 subgroups of G}
From Theorem \ref{MSSV} we  know that $[G$:$H]=2$ except for $III-c)$. Such a pair is given by an exact sequence
$$1 \to H \to G \to \IZ/2 \to 1.$$
This type of extensions, where $H=D_n$ and $G$ has another subgroup $H'$ isomorphic to $D_n$, has been classified in \cite{CLP2}, Proposition 7.4. There are 3 cases, which we call \emph{group types}:\\

\noindent
\underline{$Group\ type\ 1)$} $G\cong D_n\times \mathbb{Z}/2$, $H$ corresponds to the subgroup $D_n\times \{0\}$.\\
\underline{$Group\ type\ 2)$} $n=2d$, $G\cong D_{2n}=<z,y|z^{2n}=y^2=1,yzy=z^{-1}>$, $H=<x:=z^2,y>$.\\
\underline{$Group\ type\ 3)$} $n=4h$, where $h$ is odd, and G is the semidirect product of $H\cong D_n$ with $<\beta_2>\cong \mathbb{Z}/2$, such that conjugation by $\beta_2$ acts as follows:
$$y\mapsto yx^2,x\mapsto x^{2h-1}.$$
For each group type, we will determine the index 2 subgroups of $G$ and find out which of them are isomorphic to $D_n$.\\

\noindent
\underline{$Group\ type\ 1)$} Recall the standard presentation $D_n=<x,y|x^n=y^2=1,yxy^{-1}=x^{-1}>$ and let $C_n:=\mathbb{Z}/n$.\\
We have to understand the index 2 subgroups $K$ of $D_n$, such that $K\triangleleft G$, where $K$ corresponds to $H \cap H'$. \\
$a)$ $K= C_n\times 0$ (This is the only case when n is odd).\\
Since $G/K\cong (\mathbb{Z}/2)^2$, there are two more index 2 subgroups $H_{1,1}:=<K,(e,1)>$,\\
$H_{1,2}:=<K,(y,1)>\cong D_n$.\\
$b)$ If $n=2d$, there are two more cases, $K=<(x^2,0),(y,0)>$ or $K=<(x^2,0),(yx,0)>$ 
(both isomorphic to $D_d$). \\
Here we have 4 more index 2 subgroups, $H_{1,3}:=<(x^2,0),(y,0),(e,1)>$, $H_{1,4}:=<(x^2,0),(y,0),(x,1)>$, $H_{1,5}:=<(x^2,0),(yx,0),(e,1)>$,  $H_{1,6}:=<(x^2,0),(yx,0),(x,1)>$. On checks easily that $H_{1,4}$ and $H_{1,6}$ are isomorphic to $D_n$ and that $H_{1,3}$ and $H_{1,5}$ are isomorphic to $D_n$ if and only if $d$ is odd.\\
\underline{$Group\ type\ 2)$} Using similar arguments as for group type 1), we obtain 2 more index 2 subgroups: $H_{2,1}=C_{2n}$, $H_{2,2}=<z^2,yz> \cong D_n$.\\
\underline{$Group\ type\ 3)$}
There are 6 more index 2 subgroups: $H_{3,1}=<C_n,(e,\beta_2)>$, $H_{3,2}=<C_n,(y,\beta_2)>$, $H_{3,3}=<(x^2,0),(y,0),(e,\beta_2)>$, $H_{3,4}=<(x^2,0),(y,0),(x,\beta_2)>$, $H_{3,5}=<(x^2,0),(yx,0),(e,\beta_2)>$, $H_{3,6}=<(x^2,0),(yx,0),(x,\beta_2)>$,
 and only $H_{3,3}$ is isomorphic to $D_n$ ( since $H_{3,3}=<(y,\beta_2),(e,\beta_2)>$).

\section{Hurwitz vectors for $C \to C/G$}\label{admissible f's}
We start by recalling some general theory of Galois covers of Riemann surfaces (cf. \cite{Cat2}, section 5).\\
Let $H$ be a finite group (not necessarily isomorphic to $D_n$) which acts effectively on a curve $C$ of genus $g\geq 2$, we obtain a Galois cover $p: C\rightarrow C/H:=C'$ branched on $r$ points with branching indices $m_1,...,m_r$. Denote by $g'$ the genus of $C'$, the {\em orbifold fundamental group} of the cover is a group with the following presentation:
$$T(g';m_1,...,m_r):=<\alpha_1,\beta_1,...,\alpha_{g'},\beta_{g'};\gamma_1,...,\gamma_r\ |\ \Pi[\alpha_j,\beta_j]\cdot \Pi\gamma_i=1, \gamma_i^{m_i}=1>$$
The cover $C\rightarrow C/H$ is (topologically) determined by a surjective morphsim 
$$f: T(g';m_1,...,m_2)\rightarrow G,$$ such that $f(\gamma_j)$ has order $m_j$ inside $G$. We call $v:=[f(\alpha_1),f(\beta_1),...,f(\alpha_{g'}),f(\beta_{g'});f(\gamma_1),...,f(\gamma_r)]$ the {\em Hurwitz vector} associated to $f$.\\
In this section we study the Hurwitz vectors of each cover type $C\rightarrow C/G$ in Theorem $\ref{MSSV}$. Hence we have that $C/G\simeq \mathbb{P}^1$, and we set $T(m_1,...,m_r):=T(0;m_1,...,m_r)$. \\
Given a morphism $f:T(m_1,...,m_r)\rightarrow G$, the Hurwitz vector associated to $f$ is not uniquely determined, since we can choose different presentations for $T(m_1,...,m_r)$. For instance consider $T(m_1,...,m_r)$ with the presentation $<\gamma_1,...\gamma_r|\Pi\gamma_i=1,\gamma_i^{m_i}=1>$, for any $1\leq k<r$, we have a set of generators $\{\delta_i\}$, where $\delta_i:=\alpha_i$ if $i\neq k,k+1$; $\delta_k:=\alpha_{k}\alpha_{k+1}\alpha^{-1}_{k}$ and $\delta_{k+1}:=\alpha_k$, this induces an isomorphism between $T(m_1,...,m_r)$ and $T(l_1,...,l_r)$, where $l_i=m_i$ if $i\neq k,k+1$; $l_k=m_{k+1}$ and $l_{k+1}=m_k$. Different choices of the generators correspond to the following braid group action on the set of Hurwitz vectors.\\
Recall that Artin's {\em braid group on $r$ strands} has the presentation 
$$\mathcal B_r := <\sigma_1,...,\sigma_{r-1} | \forall 1\leq i\leq r-2,\ \sigma_i\sigma_{i+1}\sigma_i = \sigma_{i+1}\sigma_i\sigma_{i+1}; \forall   |j -i| \geq 2,\sigma_i\sigma_j=\sigma_j\sigma_i>. $$  
The group $\mathfrak{B}_r$ acts on the set of Hurwitz vectors of length $r$ as follows:
$$(v_1,...,v_i,v_{i+1},...,v_r) \stackrel{\sigma_i}{\mapsto} (v_1,...,v_iv_{i+1}v_i^{-1},v_{i},...,v_r).$$
On the other hand, for any $h\in Aut(G)$, we can compose $f$ with $h$, this induces a $Aut(G)$-action on the set of Hurwitz vectors: 
given $v=(v_1,...,v_r)$ a Hurwitz vector, define $h(v) := (h(v_1),...,h(v_r))$.\\
Since these actions (by $\mathfrak{B}_r$ and by $Aut(G)$) commute, they induce an action of the group $\mathcal B_r \times Aut(G)$ on the set of Hurwitz vectors of length $r$.
\begin{definition}\label{equivalent HV}
Given two $G$-Hurwitz vectors $v,v'$ of length $r$, we say that $v$ and $v'$ are equivalent if they are in the same $ \mathcal{B}_r\times Aut(G)$-orbit.
\end{definition}
 \begin{remark}
 Two Hurwitz vectors $v$ and $v'$ determine the same unmarked topological type iff they are equivalent (cf. \cite{CLP2}, section 2).
\end{remark}
\begin{definition}
Let $C\rightarrow C/G\cong \mathbb{P}^1$ be a Galois cover of a given group type and cover type. We call a homomorphism $f:T(m_1,\dots, m_r)\rightarrow G$ \underline{admissible} if it satisfies the following two conditions:\\
$(1)$ $f$ is surjective, $T(m_1,...,m_r)$ is isomorphic to the orbifold fundamental group of $C\rightarrow C/G$ and $f(\gamma_i)$ has order $m_i$ in $G$.\\
$(2)$ $f_{H}:=\pi_{H}\circ f: T(m_1,\dots, m_r)\rightarrow G/H$ corresponds to the cover $C/H\rightarrow \mathbb{P}^1$, where $\pi_H:G\rightarrow G/H$ is the quotient homomorphism.
\end{definition}
\begin{definition}
Let $f:T(m_1,...m_r)\rightarrow G$ and $f':T(l_1,...l_r)\rightarrow G$ be admissible for a given cover type and group type. We say $f$ is equivalent to $f'$ if their corresponding Hurwitz vectors are in the same $\mathfrak{B}_r\times Aut(G)_H$-orbit, where $Aut(G)_H$ denotes the subgroup of $Aut(G)$ which leaves $H$ invariant.
\end{definition}
\begin{remark}
An admissible $f$ determines both the covers $C\rightarrow C/G$ and $C\rightarrow C/H$, hence we require the equivalence relation to be generated by $\mathfrak{B}_r$ and $Aut(G)_H$. It can happen that two admissible homomorphisms have equivalent Hurwitz vectors, but are not equivalent (cf. Remark \ref{Rk}). 
\end{remark}

\begin{example}\label{Ex}
Cover type $III-b)$ and group type $1)$ (cf. Corollary \ref{MSSV Cor})
\end{example}
\noindent
$i)$ $c_3=2$, assume $n$ even and $c_4=n$.\\ 
Consider $f:T(2,2,2,c_4)\rightarrow D_n\times \mathbb{Z}/2$: $\gamma_1\mapsto (yx,1)$, $\gamma_2\mapsto (e,1)$, $\gamma_3\mapsto (y,0)$, $\gamma_4\mapsto (x,0)$.\\
$\delta_{H_{1,2}}=\delta_{H_{1,6}}=1$, $\delta_{H_{1,4}}=2$.\\
$ii)$ $c_3>2$, assume we have an admissible $f$, it is easy to see that $f(\gamma_3)=(x^{i_3},0)$, $f(\gamma_4)=(x^{i_4},0)$. $f(\gamma_1),f(\gamma_2)\in\{(yx^k,1),k\in \mathbb{Z}; (x^{n/2},1)(if\ n\ is\ even) \}$. Since $\Pi f(\gamma_i)=1$, there are only two possibilities:\\
$(a)$ $f(\gamma_1),f(\gamma_2)=(x^{n/2},1)$, which implies $\mbox{Im}(f)\subset <(x,0),(0,1)>$, a contradiction.\\
$(b)$ $f(\gamma_1)=(yx^{i_1},1)$, $f(\gamma_2)=(yx^{i_2},1)$, which implies $\mbox{Im}(f)\subset <(x,0),(y,1)>$, again a contradiction.\\

Now we classify all admissible $f$'s  for the covering $C \to C/G$,  in the following way: For each cover type and group type, we construct all possible Hurwitz vectors according to their branching behavior, as given in Theorem \ref{MSSV}. \\

\begin{lemma}
Group type $2)$ has no admissible $f$ for any cover type.
\begin{proof}
 \underline{Cover type I)} \\
Assume we have an admissible $f:T(2,2,2,2,2,2)\rightarrow D_{2n}$, then $f_{H}(\gamma_i)=1$, $i=1,\dots,6$, which implies that
 $f(\gamma_i)\in \{yz^{2k+1}, z^{2l+1},k,l\in \mathbb{Z} \}$. Moreover $f(\gamma_i)$ has order two, thus $f(\gamma_i)\in\{yz^{2k+1},k\in \mathbb{Z} \}$. We find that $\mbox{Im}(f)\subset H_{2,2}$, a contradiction.\\
 \underline{ Cover type II) }\\
If there exists an admissible $f:T(2,2,2,2,c_5)\rightarrow D_{2n}$, we get $f(\gamma_i)\in \{yz^{2k+1},k\in \mathbb{Z} \}, i=1,2,3,4$ and $f(\gamma_5)\in \{z^{2l},l\in \mathbb{Z} \} $(since $\Pi f(\gamma_i)=1$), which implies that $\mbox{Im}(f)\subset H_{2,2}$, a contradiction.\\
 \underline{Cover type III-a)} \\
Given an admissible $f:T(2,2,2,2d_4)\rightarrow D_{2n}$, we get $f(\gamma_i)\in\{yz^{2k+1},k\in \mathbb{Z} \}, i=1,2,3$, and $f(\gamma_4)\in \{z^{2l+1},l\in \mathbb{Z} \}$. However, $\Pi f(\gamma_i)\neq 1$, a contradiction.\\
 \underline{Cover type III-b)}\\
$i)$ $c_3=2$. We have $f(\gamma_i)=yz^{2k_i+1},i=1,2, f(\gamma_3)=yz^{2k_3}$ or $z^n$,  $f(\gamma_4)=z^{2k_4}$. If $f(\gamma_3)=yz^{2k_3}$ we find $\Pi f(\gamma_i)\neq 1$; otherwise $f(\gamma_3)=z^n$, which implies $\mbox{Im}(f)\subset <yz,z^2>$. In both cases we have no admissible $f$.\\
$ii)$ $c_3>2$. We have $(f(\gamma_1),f(\gamma_2),f(\gamma_3),f(\gamma_4))=(yz^{2k_1+1},yz^{2k_2+1},z^{2k_3},z^{2k_4})$. We see $\mbox{Im}(f)\subset <yz,z^2>$, a contradiction.
\end{proof}
\end{lemma}

\begin{lemma}\label{group type 3}
Group type $3)$ has no admissible $f$ for any cover type.
\begin{proof}
  First we determine the order 2 elements of type $(a,\beta_2)$ in $G$. One computes easily that $(x^j,\beta_2)^2=(x^{2jh},0)$ and $(yx^k,\beta_2)^2=(x^{2kh-2k+2},0)\neq (e,0)$. Therefore we conclude that $(a,\beta_2)$ is of order two $\Leftrightarrow$ $a=x^j$ and $j$ is even. \\
 \underline{Cover type I)}\\ 
Now assume we have an admissible $f$, which implies that $f(\gamma_i)=(x^{2j_i},\beta_2)$. However these elements are contained in the proper subgroup $<(x^{2},0),(e,\beta_2)>$, we see $f$ can not be surjective, a contradiction.\\
 \underline{Cover type II)}\\
If there exists an admissible $f$, we must have $f(\gamma_i)=(x^{2j_i},\beta_2),i=1,2,3,4$, and since $\Pi f(\gamma_i)=1$ it follows that $\mbox{Im}(f)\subset <(x^2,0),(e,\beta_2)>$, a contradiction. \\
 \underline{Cover type III-a)}\\
Assume we have an admissible $f$, we see that $f(\gamma_i)=(x^{2j_i},\beta_2),i=1,2,3$. Since $\Pi f(\gamma_i)=1$ it follows that $\mbox{Im}(f)\subset <(x^2,0),(e,\beta_2)>$, again a contradiction. \\
 \underline{Cover type III-b)}\\ 
$i)$ $c_3=2$.   We must have $f(\gamma_1)=(x^{2j_1},\beta_2)$, $f(\gamma_2)=(x^{2j_2},\beta_2)$, $f(\gamma_3)=(x^{2h},0)$ or $(yx^k,0)$, $f(\gamma_4)=(x^l,0),l\neq 2h$. If $f(\gamma_3)=(x^{2h},0)$, then $\mbox{Im}(f)\subset <(x,0),(0,\beta_2)>$; if $f(\gamma_3)=(yx^k,0)$ we see $\Pi f(\gamma_i)\neq 1$. In both cases we can not get an admissible $f$.\\
$ii)$ $c_3>2$. Given an admissible $f$, we have $f(\gamma_1)=(x^{2j_1},\beta_2)$, $f(\gamma_2)=(x^{2j_2},\beta_2)$, $f(\gamma_3)=(x^{k_3},0)$ and
$f(\gamma_4)=(x^{k_4},0)\ (k_3,k_4\neq 2h)$. One sees immediately that $\mbox{Im}(f)\subset <(x,0),(0,\beta_2)>$, a contradiction.
\end{proof}
\end{lemma} 

\begin{lemma}
Cover type $III-c)$ has no admissible f.
\begin{proof}
Assume that we have an admissible $f:T(2,2,2,c_4)\rightarrow G$. \\
Let $(b_1,b_2,b_3,b_4):=(f(\gamma_1),f(\gamma_2),f(\gamma_3),f(\gamma_4))$. We have \\
$(1)$ $b_1^2=b_2^2=b_3^2=1$. Since $b_4\in H$ and $order(b_4)=c_4>2$, we see that $b_4$ must lie in the cyclic group,  say $b_4=x^k$, we also find $n>2$.\\
$(2)$ The fact that $H$ is normal in $G$ implies that $b_ixb_i=x^{k_i},i=1,2,3$, therefore $x^kb_i=b_ix^{kk_i},i=1,2,3$. \\
$(3)$ $b_1b_2b_3b_4=1\Rightarrow$ $b_1b_2=x^{-k}b_3$, moreover $(b_1b_2)^2=x^{-k}b_3x^{-k}b_3=x^{-k-kk_3}$. \\
Any element in $\mbox{Im}(f)$ has the form $\Pi_{i=1}^l \beta_i $, where $\beta_i\in \{x^k,x^{-k}, b_1,b_2,b_3\}$. Since $b_1b_2b_3b_4=1$, without loss of generality we can assume $\beta_i\in \{x^k,x^{-k}, b_1,b_2\}$, which means that every element in $\mbox{Im}(f)$ is a word in these four elements. \\
Using (2), we can "move" the $x^{\pm k}$ terms to the end. Taking $(1)$ into account, we see that the elements are of the forms $(b_1b_2)^sx^t$, $b_2(b_1b_2)^sx^t$ or
 $(b_1b_2)^sb_1x^t$, now use $(3)$, one sees immediately that elements in $\mbox{Im}(f)$ have the form $x^j$, $b_1x^j$,  $b_2x^j$ or $b_3x^j$. It turns out that $H\nsubset \mbox{Im}(f)$, a contradiction.
\end{proof}
\end{lemma}

From the preceeding, we know that the only group type to consider is Group type I). We denote by $(e,0)$ the neutral element of $D_n \times \IZ/2$, where $\IZ/2$ is additively generated by $1$.\\

For the action of the braid group on the set of Hurwitz vectors we make use of Lemma 2.1 in \cite{CLP1}.

\begin{lemma}\label{Lemma 2.1}

Every Hurwitz vector of length $r$  with elements in $D_n$ of the form $$v=(v_1,...,yx^a,yx^b,yx^c,...,v_r)$$ is equivalent to $v'=(v_1,...,yx^{a'},yx^{a'},yx^{c'},...,v_r)$ or $v''=(v_1,...,yx^{a'},yx^{b'},yx^{b'},...,v_r)$ via braid moves that only affect the triple $(yx^a,yx^b,yx^c)$. 

\end{lemma}

\begin{lemma} Classification of cover type I)\

In this case the only admissible Hurwitz vector for \textbf{$n$ odd} is

\[ v = ((y,1),(y,1),(yx,1),(yx,1),(e,1),(e,1)).\]
 
For \textbf{$n$ even} (n=2m) there are the following possibilities: 

\[v = ((y,1),(y,1),(yx,1),(yx,1),(e,1),(e,1)),\]


\[v = ((y,1),(yx^{m},1),(yx,1),(yx,1),(x^{m},1),(e,1)),\]

\[v = ((y,1),(yx^{m},1),(yx^2,1),(yx^2,1),(x^{m},1),(e,1)), m ~\mbox{odd.}\]

For $n=2$ there are the following:

\[ v = ((y,1),(y,1),(x,1),(x,1),(e,1),(e,1)),\]
\[v = ((y,1),(yx,1),(x,1),(x,1),(x,1),(e,1))\]
\begin{proof}

Since the cover $C/H \to \IP^1$ branches in 6 points (cf. [MSSV]) we need a Hurwitz vector with second component equal to $1$. So we have 

\[v = ((y^{k_1}x^{l_1},1),(y^{k_2}x^{l_2},1),(y^{k_3}x^{l_3},1),(y^{k_4}x^{l_4},1),(y^{k_5}x^{l_5},1),(y^{k_6}x^{l_6},1))\] 

The first observation is that the condition $<v> = G$ implies that there must exist $j$, s.t. $k_j = 1$.
Therefore up to automorphism we can assume 

\[v = ((y,1),(y^{k_2}x^{l_2},1),(y^{k_3}x^{l_3},1),(y^{k_4}x^{l_4},1),(y^{k_5}x^{l_5},1),(y^{k_6}x^{l_6},1))\] 

We consider the two cases $n$ odd and $n$ even separately.

\begin{itemize}

\item[i)] \underline{n odd}: Not all $k_j$ can be equal to $1$. Otherwise we cannot generate the element $(y,0)$. Now the only element of order two of the form $(x^l,1)$ in $G$ is $(e,1)$. So because of the product one condition $v$ either looks like

\[v = ((y,1),(yx^{l_2},1),(yx^{l_3},1),(yx^{l_4},1),(e,1),(e,1))\] 

or

\[v = ((y,1),(y,1),(e,1),(e,1),(e,1),(e,1)),\] 

the latter being excluded, since $G\neq <v>$.\\

The product one condition gives $l_2 + l_4 \equiv l_3 \mod n$. The condition $<v> = G$ implies $\gcd(l_2,l_3,l_4,n) = \gcd(l_2,l_4,n) = 1$. Since the second factor $\IZ/2$ of $G$ is abelian, we can apply Lemma \ref{Lemma 2.1} to achieve that $l_3 = l_4$. Now $v$ looks like 

\[ v = ((y,1),(yx^{l_2},1),(yx^{l_4},1),(yx^{l_4},1),(e,1),(e,1)) \]

and again by product one we obtain $l_2 \equiv 0 \mod n$ and therefore $1 = \gcd(l_2,l_4,n) = \gcd(l_4,n)$. \\

So we can apply the automorphism $(x^{l_4},0)  \mapsto (x,0) ,(y,0) \mapsto (y,0)$ to $v$ and we can take 

\[ v = ((y,1),(y,1),(yx,1),(yx,1),(e,1),(e,1)) \]

as a Hurwitz vector for the covering $C \to \IP^1$.

\item[ii)]\underline{n even}: Recall the general form:

\[v = ((y,1),(y^{k_2}x^{l_2},1),(y^{k_3}x^{l_3},1),(y^{k_4}x^{l_4},1),(y^{k_5}x^{l_5},1),(y^{k_6}x^{l_6},1))\]

Again, first we distinguish the possible Hurwitz vectors by the (even and positive) number of $k_{j}$ that are equal to $1$. We call the element $y^kx^l$ a reflection if $k\equiv 1(\mod2).$\\
In the current case there exists $m = n/2$, which gives the extra order $2$ element $(x^m,1) \in G$. As in the odd case, 6 reflections cannot occur. For the case of 2 reflections, assume, up to ordering,

\[v = ((y,1),(yx^{l_2},1),(x^{l_3},1),(x^{l_4},1),(x^{l_5},1),(x^{l_6},1)).\]

As before, $(l_3,l_4,l_5,l_6) = (0,0,0,0)$ is impossible. In the cases $(l_3,l_4,l_5,l_6) = (m,m,0,0)$ and $(l_3,l_4,l_5,l_6) = (m,m,m,m)$ we get $l_2 = 0$. In the first case we can only have $<v> = G$ if $n=2$. Also in the second case we must have $n=2$ but the elements $(y,1)$ and $(x,1)$ cannot generate $G$ since the element $(e,1)$ is missing. In the cases $(l_3,l_4,l_5,l_6) = (m,m,m,0)$ and $(l_3,l_4,l_5,l_6) = (m,0,0,0)$ we get $l_2 = m$ which also implies that $n=2$. So if $n>2$ these cases don't occur. The corresponding Hurwitz vectors are:

\[v = ((y,1),(y,1),(x,1),(x,1),(e,1),(e,1)),\]
\[v = ((y,1),(yx,1),(x,1),(x,1),(x,1),(e,1))\]
and
\[v = ((y,1),(yx,1),(x,1),(e,1),(e,1),(e,1)),\]

the third one being equivalent to the second one by an automorphism of $G$ that fixes $D_n$.

Assume, for the case of 4 reflections, up to ordering

\[v = ((y,1),(yx^{l_2},1),(yx^{l_3},1),(yx^{l_4},1),(x^{l_5},1),(x^{l_6},1)).\]

Here we have the 3 cases: $l_5 = l_6 = m$, $l_5 = l_6 = 0$ and $l_5 = m$, $l_6 = 0$.\\

In the first 2 cases from the product-one condition we get $l_2 + l_4 \equiv l_3 \mod n$.
To generate $G$ we must have $\gcd(l_2,l_3,l_4,n) = \gcd(l_2,l_4,n) = 1$.

Using Lemma \ref{equivalent HV} again, we arrive at 

\[v = ((y,1),(yx^{l_2},1),(yx^{l_4},1),(yx^{l_4},1),(x^{m},1),(x^{m},1))\]

resp.

\[v = ((y,1),(yx^{l_2},1),(yx^{l_4},1),(yx^{l_4},1),(e,1),(e,1))\]

and so we get $l_2 \equiv 0 \mod n$. Now we have $\gcd(l_2,l_4,n) = \gcd(l_4,n) = 1$ and we can apply the automorphism $x^{l_4}  \mapsto x, y \mapsto y$ to $v$ to arrive at

\[v = ((y,1),(y,1),(yx,1),(yx,1),(x^{m},1),(x^{m},1))\]

resp.

\[v = ((y,1),(y,1),(yx,1),(yx,1),(e,1),(e,1)).\]

Using the morphism $(e,1) \mapsto (x^m,1), (y,0) \mapsto (yx^{-m},0)$ we see that these two are equivalent.

It remains to consider the case $l_5 = m$ and $l_6 = 0$, i.e.

 \[v = ((y,1),(yx^{l_2},1),(yx^{l_3},1),(yx^{l_4},1),(x^{m},1),(e,1)).\]
 
 We apply Lemma 2.1, [CLP1] again and it follows $l_2 = m$. So we get
 
\[v = ((y,1),(yx^{m},1),(yx^{l},1),(yx^{l},1),(x^{m},1),(e,1))\]

where $\gcd(l,m) = 1$.\\

We have two sub cases, i.e. $\gcd(l,n) = 1$ and $\gcd(l,n) = 2$. In the first case we can use the automorphism $x^l \mapsto x, y \mapsto y$ to obtain 
\[v = ((y,1),(yx^{m},1),(yx,1),(yx,1),(x^{m},1),(e,1)).\]

In the second case (where $m$ must be odd) we can achieve

\[v = ((y,1),(yx^{m},1),(yx^2,1),(yx^2,1),(x^{m},1),(e,1)).\]

\end{itemize} 

\end{proof}

\end{lemma}

\begin{lemma}
Classification of cover type II)\\
Up to equivalence, the admissible $f$ is given by the Hurwitz vector:\\
$(1)$ $c_5=2$, 
$$v=((y,1),(yx,1),(yx,1),(e,1),(y,0)),$$
$(2)$ $c_5>2$,
$$v=((y,1),(yx^{-1},1),(e,1),(e,1),(x,0)),c_5=n,$$
$$v=((y,1),(yx^{m-1},1),(x^m,1),(e,1),(x,0)),n=2m,c_5=n,$$
$$v=((y,1),(yx^{m-2},1),(x^m,1),(e,1),(x^2,0)),n=2m,m\ is\ odd,c_5=m,$$
\begin{proof}
 Assume we have an admissible $f:T(2,2,2,2,c_5)\rightarrow D_n\times \mathbb{Z}/2$.\\
we must have:
$$v:=(f(\gamma_1),f(\gamma_2),f(\gamma_3),f(\gamma_4),f(\gamma_5))=((a_1,1),(a_2,1),(a_3,1),(a_4,1),(a_5,0))$$
There are two cases:\\
$(1)$ $\underline{c_5=2}$.\\
As in the previous argument, we do the classification in terms of the number of reflections in $\{a_i\}$, which can be either 2 or 4.\\
$(i)$ There are 2 reflections.\\
$(a)$ $a_5$ is a reflection, W.L.O.G we can assume $a_1$ is another reflection, and $a_1=yx^l,a_5=y$. $a_2,a_3,a_4\in \{e,x^{n/2}$(if n is even)$\}$.\\
There are 4 cases (up to an order change): $\alpha)\ (a_2,a_3,a_4)=(e,e,e)$, $\beta)\ (a_2,a_3,a_4)=(x^{n/2},e,e)$, $\gamma)\ (a_2,a_3,a_4)=(x^{n/2},x^{n/2},e)$, $\delta)\ (a_2,a_3,a_4)=(x^{n/2},x^{n/2},x^{n/2})$.\\
Case $\alpha)$,$\delta)$ we get no admissible $f$ since $f$ can not be surjective.\\
For case $\beta),\gamma)$ (where n is even) we get $f$ is admissible $\iff$ $n=2$.\\
$(b)$ $a_5$ is not a reflection, first we conclude that $n$ must be even and $a_5=x^{n/2}$. Using similar arguments as in $a)$, one finds that
$$v=((y,1),(yx^l,1),(a_3,1),(a_4,1),(x^{n/2},0)),\ a_3,a_4\in \{e,x^{n/2}\}.$$
There are three cases, and one checks easily that in each case $f$ is admissible if and only if $n=2$.\\
$(ii)$ There are 4 reflection.\\
$a)$ $a_5$ is a reflection. W.L.O.G we assume
$$v=((yx^{l_1},1),(yx^{l_2},1),(yx^{l_3},1),(a_4,1),(y,0)),\ a_4\in \{e, x^{n/2}\textit{(if n is even)}\}.$$
Again we apply Lemma \ref{equivalent HV} so that we can assume $l_2=l_3$. Since $f$ is admissible, (using similar arguments as in the previous Lemma,) 
we have:\\
Case $\alpha)$ If $a_4=e$, then $l_1\equiv 0\ (n)$, $\gcd(l_2,n)=1$. Under the automorphism $x^{l_2}\mapsto x$, $y\mapsto y$, we get
$$v\sim((y,1),(yx,1),(yx,1),(e,1),(y,0)).$$
Case $\beta)$ $n=2m$ and $a_4=x^m$. One gets $l_1\equiv m\ (2m)$, and $\gcd(l_2-m,2m)=1$. Using the automorphism $x^{l_2-m}\mapsto x$, $y\mapsto y$, then we can achieve\\
$$v\sim((yx^m,1),(yx^{m+1},1),(yx^{m+1},1),(x^m,1),(y,0)).$$
Using the automorphism (of $G$): $(x,0)\mapsto (x,0)$, $(y,0)\mapsto (y,0)$, $(e,1)\mapsto (x^m,1)$, one finds that Case $\beta)$ is equivalent to Case $\alpha)$.\\
$b)$ $a_5$ is not a reflection.\\
In this case $n$ must be even, and $v=((y,1),(yx^{l_2},1),(yx^{l_3},1),(yx^{l_4},1),(x^{n/2},0))$. It is easy to see that $f$ can not be surjective since $(y,0)$ is not contained in the image.\\
Up to now we have got all the admissible $f$'s for the case $n=2$. (Since $n=2$ implies that $c_5=2$). One checks easily that they are equivalent to each other, since in this case $G$ is abelian.\\
$(2)$ $\underline{c_5>2}$.\\
$a_5$ must lie in the cyclic subgroup, say $a_5=x^k\ (k\neq \frac{n}{2}$ if n is even $)$.  \\
$(i)$ There are 2 reflections, W.L.O.G. we assume
$$v=((y,1),(yx^l,1),(a_3,1),(a_4,1),(x^k,0)),\ a_3,a_4\in \{e,x^{n/2}(if\ n\ is\ even)\}$$
There are 3 cases:\\
Case $\alpha)$ $(a_3,a_4)=(e,e)$.\\
We get $l+k\equiv 0\ (\mod n)$ and $\gcd(k,n)=1$. Applying the automorphism $x^k\mapsto x$, $y\mapsto y$ we get
$$v\sim((y,1),(yx^{-1},1),(e,1),(e,1),(x,0)).$$
Moreover we see that $c_5=n$.\\
Case $\beta)$ $n=2m$ and $(a_3,a_4)=(x^m,e)$.\\
We get $l+k\equiv m\ (\mod 2m)$ and $\gcd(k,m)=1$. \\
If $\gcd(k,n)=1$ (which is the unique case if $2|m$), 
$$v\sim ((y,1),(yx^{m-1},1),(x^m,1),(e,1),(x,0))$$
Here we find $c_5=n$.\\
Otherwise $\gcd(k,n)=2$ (which may happen only when $2\nmid m$),
$$v\sim ((y,1),(yx^{m-2},1),(x^m,1),(e,1),(x^2,0))$$
and we have $c_5=m$.\\
Case $\gamma)$ $n=2m$ and $(a_3,a_4)=(x^m,x^m)$.\\
We get $l+k\equiv 0\ (n)$ and $\gcd(k,n)=1$.
$$v\sim ((y,1),(yx^{-1},1),(x^m,1),(x^m,1),(x,0)),c_5=n$$
Using the automorphism $(x,0)\mapsto (x,0)$, $(y,0)\mapsto (yx^{-m},0)$, $(e,1)\mapsto (x^m,1)$, one finds case $\gamma)$ is equivalent to Case $\alpha)$.\\
$(ii)$ There are 4 reflections.\\
One checks easily that $f$ can not be surjective since $(y,0)\notin \mbox{Im}(f)$.
\end{proof}
\end{lemma}

\begin{lemma}
Classification of cover type III-a)\\
We have that $n=2m$ and $d_4=m$. Up to equivalence there is a unique admissible $f$ given by the Hurwitz vector:
$$v=((y,1),(yx^{-1},1),(e,1),(x,1)).$$
\begin{proof}
Assume $f:T(2,2,2,2d_4)\rightarrow D_n\times \mathbb{Z}/2$ is admissible.
$$v:=(f(\gamma_1),f(\gamma_2),f(\gamma_3),f(\gamma_4))=((a_1,1),(a_2,1),(a_3,1),(a_4,1)).$$
$d_4>1\Rightarrow a_5=x^k$ $(k\neq n/2\ if\ n\ is \ even)$.\\
There can only be 2 reflections among $a_1,a_2,a_3$. W.L.O.G. we can assume
$$v=((y,1),(yx^l,1),(a_3,1),(x^k,1)),a_3\in \{e,n/2(if\ n\ is\ even)\}$$
Case $a)$ $a_3=e$.\\
We get $l+k\equiv 0\ (\mod n)$ and $\gcd(k,n)=1$,
$$v\sim ((y,1),(yx^{-1},1),(e,1),(x,1))$$
In this case $2d_4=n$, it turns out that $n$ must be even.\\
Case $b)$ $n=2m$ and $a_3=x^{m}$.\\
We get $l+k\equiv m\ (2m)$ and $\gcd(l,n)=1$,
$$v\sim ((y,1),(yx,1),(x^m,1),(x^{m-1},1)).$$
Using the automorphism $(x,0)\mapsto (x^{-1},0)$, $(y,0)\mapsto (yx^{-m},0)$, $(e,1)\mapsto (x^m,1)$, we find that Case $b)$ is equivalent to Case $a)$.
\end{proof}
\end{lemma}

\begin{lemma}
Classification of cover type III-b)\\
We have that $c_3=2$ and $c_4=n$. Up to equivalence there is a unique  admissible $f$ given by the Hurwitz vector:
$$v=((yx,1),(e,1),(y,0),(x,0)).$$
\begin{proof}
From Example \ref{Ex} we see if that a type $III-b)$ cover has group type $1)$, $c_3$ must be 2, combining with the proof of Corollary  \ref{MSSV Cor} one obtains that the case $(\delta_H,\delta_{H'})=(1,3)$ does not occur. \\
Let $f:T(2,2,2,c_4)\rightarrow D_n\times \mathbb{Z}/2$ be admissible. We must have
$$v:=(f(\gamma_1),f(\gamma_2),f(\gamma_3),f(\gamma_4))=((a_1,1),(a_2,1),(a_3,0),(a_4,0))$$
Since $c_4>2$ we get $a_4=x^k$. It is obvious that there are two (and only two) reflections among $a_1,a_2,a_3$.\\
$(1)$ $a_3$ is not a reflection. $n$ must be even (let $n=2m$) and $a_3=x^m$.
W.L.O.G we assume 
$$v=((y,1),(yx^l,1),(x^m,0),(x^k,0)).$$
It is easy to see $(y,0)\notin \mbox{Im}(f)$, therefore in this case there is no admissible $f$.\\
$(2)$ $a_3$ is a reflection. W.L.O.G we assume 
$$v=((yx^l,1),(a_2,1),(y,0),(x^k,0)),a_2\in \{e,n/2(if\ n\ is\ even)\}.$$
$(i)$ $a_2=e$, we get $k\equiv l\ (n)$ and $\gcd(k,n)=1$,
$$v\sim (yx,1),(e,1),(y,0),(x,0)), c_4=n.$$
$(ii)$ $n=2m$ and $a_2=x^m$, we get $k\equiv l+m\ (2m)$, $\gcd(k,n)=1$,
$$v\sim (yx^{m+1},1),(x^m,1),(y,0),(x,0)), c_4=n.$$
Using the automorphism $(x,0)\mapsto (x,0)$, $(y,0)\mapsto (y,0)$, $(e,1)\mapsto (x^m,1)$, we see that Case $(ii)$ is equivalent to Case $(i)$.
\end{proof}
\end{lemma}

\begin{remark} \label{Rk}
If we drop the restriction on $f_H$, it is easy to check that the Hurwitz vectors in $III-a)$ and $III-b)$ are equivalent. 
(Consider the automorphism of $G$: $(x,0)\mapsto (x,1)$, $(y,0)\mapsto (yx,0)$, $(e,1)\mapsto (e,1)$)
\end{remark}
\section{Acknowledgement}
The present work took mainly place in the realm of the DFG Forschergruppe 790 "Classification of algebraic surfaces and compact complex manifolds". 
The first author is currently sponsored by the project "ERC Advanced Grant 340258 TADMICAMT".\\
We would like to thank Fabrizio Catanese for suggesting the topic of this paper. We also would like to thank Michael L{\"o}nne and Fabrizio Catanese for carefully reading the paper and several useful suggestions.

\bibliographystyle{alpha}

\noindent
Authors' Address: 

\noindent
Binru Li, Sascha Weigl\\
Lehrstuhl Mathematik VIII, Universit\"at Bayreuth\\
Universit\"atsstra\ss e  30, D-95447 Bayreuth\\
E-mail address:\\
binru.li@uni-bayreuth.de\\
sascha.weigl@uni-bayreuth.de

\end{document}